\documentclass[oneside,12pt,letterpaper]{amsart}

\usepackage{amsmath,amssymb,epsf}
\usepackage{amsfonts,amsbsy,amscd,stmaryrd}
\usepackage{latexsym,amsopn}
\usepackage{amsthm}
\usepackage{mathrsfs}
\usepackage[pagebackref=false,bookmarks=false,pdfborder={0 0 0}]{hyperref}

\allowdisplaybreaks
 
\usepackage{graphicx}

\usepackage{color}
\definecolor{Red}{rgb}{1.,0.,0.}
\newcounter{smallarabics}
\newenvironment{arabicenumerate}
{\begin{list}{{\normalfont\textrm{(\arabic{smallarabics})}}}
  {\usecounter{smallarabics}\setlength{\itemindent}{0cm}
   \setlength{\leftmargin}{5ex}\setlength{\labelwidth}{4ex}
   \setlength{\topsep}{0.75\parsep}\setlength{\partopsep}{0ex}
   \setlength{\itemsep}{0ex}}}
{\end{list}}

\newcounter{smallroman}

\newcommand{\ben}{\begin{arabicenumerate}}  
\newcommand{\een}{\end{arabicenumerate}}

\newtheorem{theorem}{Theorem}[section]

\newtheorem{proposition}[theorem]{Proposition}
\newtheorem{lemma}[theorem]{Lemma}

\theoremstyle{definition}
\newtheorem{definition}[theorem]{Definition}
\newtheorem{remark}[theorem]{Remark}
\newtheorem{example}[theorem]{Example}
\newcommand{\beq}{\begin{equation}}
\newcommand{\eeq}{\end{equation}}

\newcommand{\bea}{\begin{aligned}}
\newcommand{\eea}{\end{aligned}}

\newcommand{\bex}{\begin{example}}
\newcommand{\eex}{\end{example}}
\def\bel{\begin{lemma}}
\def\eel{\end{lemma}}
\def\bet{\begin{theoreme}}
\def\eet{\end{theoreme}}
\def\bed{\begin{definition}}
\def\eed{\end{definition}}
\def\ber{\begin{remark}}
\def\eer{\end{remark}}

\def\rr{{\mathbb R}}

\def\nn{{\mathbb N}}

\def\part{{\rm par}}

\usepackage{calrsfs}
\DeclareMathAlphabet{\pazocal}{OMS}{zplm}{m}{n}

\def\cY{{\pazocal Y}}

\def\cN{{\pazocal N}}

\def\wf{{\rm WF}}

\DeclareMathOperator{\Ker}{Ker}

\def \p{ \partial}

\def\12{\frac{1}{2}}
\def\14{\frac{1}{4}}

\DeclareMathOperator{\supp}{supp}

\newcommand{\one}{\boldsymbol{1}}

\def\cX{{\pazocal X}}

\def\12{\frac{1}{2}}

\def\bep{\begin{proposition}}
\def\eep{\end{proposition}}

\newcommand{\bra}{\langle} 
\newcommand{\ket}{\rangle}

\DeclareSymbolFont{boldoperators}{OT1}{cmr}{bx}{n}
\SetSymbolFont{boldoperators}{bold}{OT1}{cmr}{bx}{n}


\def\cf{C^\infty}

\def\zero{{\rm\textit{o}}}

\let\origmaketitle\maketitle
\def\maketitle{
  \begingroup
  \def\uppercasenonmath##1{} 
  \let\MakeUppercase\relax 
	\origmaketitle
  \endgroup
}

\renewcommand{\hbar}{{}}

\newcommand{\open}[1]{\mathopen{}\mathclose{\left]#1 \right[}}

\begin{document}

\title[The Feynman problem for the Klein--Gordon equation]{\large The Feynman problem for the Klein--Gordon equation}

\address{Universit\'e Paris-Saclay, D\'epartement de Math\'ematiques, 91405 Orsay Cedex, France}
\email{christian.gerard@math.u-psud.fr}
\author{\normalsize Christian G\'erard  \& Micha{\l} Wrochna}
\thanks{\emph{Acknowledgments.} Support from the grant ANR-16-CE40-0012-01 is gratefully acknowledged.}
\address{CY Cergy Paris Université, Département de Mathématiques,
2 av.~Adolphe Chauvin, 95302 Cergy-Pontoise cedex, France}
\email{michal.wrochna@cyu.fr}
\keywords{Klein--Gordon equation, scattering theory, Feynman propagators, Quantum Field Theory on curved spacetimes}
\subjclass[2010]{81T13, 81T20, 35S05, 35S35}
\begin{abstract}
We report on the well-posedness of the Feynman problem for the Klein--Gordon equation on asymptotically Minkowski spacetimes. The main result is the invertibility of the  Klein--Gordon operator with Feynman conditions at infinite times. Furthermore, the inverse is shown to coincide with the   Duistermaat--H\"ormander   Feynman parametrix modulo smoothing terms. 
\end{abstract}

\maketitle

\section{Introduction and main result}

\subsection{Introduction} For a hyperbolic partial differential operator $P$, what is the best analogue of the boundary value $(-\Delta + \lambda  + i 0)^{-1}$  of the Laplacian resolvent?

This question has a clear answer if $P=\partial_t^2-\Delta_x + m^2$ is the Klein--Gordon or wave operator on Minkowski space. Namely, if  $P_{\rm F}^{-1}$ is the Fourier multiplier by
$$
\frac{-1}{\tau^2-\xi^2-m^2 + i 0 }
$$
(where $\tau,\xi$ are the dual variables of $t,x$) then it is a formal inverse of $P$ in the sense that on test functions, $P P_{\rm F}^{-1}$ and $P_{\rm F}^{-1}  P$  are the identity. Furthermore, $P$ is easily seen to be essentially self-adjoint, and it is also not very difficult to show the \emph{limiting absorption principle}  $P_{\rm F}^{-1}=\lim_{\varepsilon\to 0^+}(P-i \varepsilon)^{-1}$ on suitable weighted $L^2$ spaces (see e.g.~\cite{derezinskisiemssen0}). The Schwartz kernel of  $P_{\rm F}^{-1}$  is called the \emph{Feynman propagator}, and plays a fundamental r\^{o}le in Quantum Field Theory. We stress that $P_{\rm F}^{-1}$ differs from the backward and forward (or advanced and retarded) formal inverses, which are the two Fourier multipliers
$$
\frac{-1}{(\tau\pm i 0 )^2-\xi^2-m^2 }
$$
and are associated to solving $Pu=f$ with $u$ or $f$ vanishing in the future, resp.~in the past (and which are therefore directly related to solving a Cauchy problem). The physical interpretation is best seen when writing the Schwartz kernel of $P_{\rm F}^{-1}$ in the time variables only: this yields the operator-valued kernel 
\beq\label{eq.feynman1}
  P_{\rm F}^{-1}(t_1,t_2)=   \frac{e^{i|t_1-t_2|\sqrt{-\Delta_x + m^2}}}{\sqrt{-\Delta_x+m^2}}.
\eeq
Even better, we can consider the Cauchy evolution of $P$ and denote by $H$ its self-adjoint generator acting on the energy space. This amounts to considering instead of $P$ the operator  
$$
\widetilde{P}=D_t- H, \quad H=\begin{pmatrix} 0 & 1 \\ -\Delta_x + m^2 & 0 \end{pmatrix},
$$
where $D_t=i^{-1}\p_t$. In this equivalent setting the Feynman propapagor is given by the formula
\beq\label{eq.feynman2}
\widetilde{P}_{\rm F}^{-1}(t_1,t_2) =\left(\one_{\rr_+}(t_1-t_2) \one_{\rr_+}(H)-\one_{\rr_-}(t_1-t_2) \one_{\rr_-}(H)\right) e^{i(t_1-t_2)H},
\eeq
where $\one_{\rr_\pm}$ is the characteristic function of the positive/negative half-line.   This elucidates Feynman's original interpretation : particles (corresponding to the positive spectral subspace of $H$) travel towards the future, whereas anti-particles (corresponding to the negative spectral subspace of $H$) travel toward the past. 

Let now $P=-\Box_g + m^2$ be the wave or Klein--Gordon operator on a Lorentzian manifold $(M,g)$. It is not difficult to imagine that the formulae \eqref{eq.feynman1} or \eqref{eq.feynman2} and the limiting absorption principle are still valid whenever the coefficients of $P$ are $t$-independent in an appropriate sense, provided that the generator $H$ has good spectral properties, see e.g.~\cite{derezinskisiemssen0} for a systematic analysis. However, outside of this exceptional case, the problem of giving meaning to the Feynman inverse $P_{\rm F}^{-1}$ of the wave or Klein--Gordon operator  $P=-\Box_g + m^2$ on a Lorentzian manifold $(M,g)$ has remained  open for a long time. A partial answer was provided by Duistermaat--H\"{o}rmander, who proposed a microlocal characterization of $P_{\rm F}^{-1}$ and who  constructed a parametrix, unique modulo smoothing terms \cite{DH}. 

Recently, the idea that emerged is that if the spacetime $(M,g)$ is not completely arbitrary, but  \emph{asymptotic} to Minkowski space, then an unambiguous definition of $P_{\rm F}^{-1}$ as a Hilbert space inverse of $P$ should be possible by imposing conditions at \emph{infinite times}. In fact, in these situations the particle and anti-particle projections $\one_{\rr^\pm}(H)$ still make sense for the asymptotic dynamics. With Feynman's interpretation in mind, if one uses the particle projection $\one_{\rr^+}(H)$ at $t=+\infty$, then the corresponding anti-particle projection $\one_{\rr^-}(H)$ should be used at $t=-\infty$. This leads to very non-local, asymptotic conditions, and an added difficulty is that even if one manages to prove invertibility on the resulting function spaces, it is not immediately clear how to check consistency with the microlocal characterisation of Duistermaat--H\"{o}rmander. In other words, one needs ways to control regularity of solutions of an inhomogeneous problem $Pu=f$ in terms of asymptotic data, and the invertibility properties are also tied to the  decay at spatial infinity. 

The first result on this kind of \emph{Feynman problem} is due to Gell-Redman--Haber--Vasy \cite{GHV}, who showed the Fredholm property of the wave operator  on asymptotically Minkowski spacetimes, later on improved by Vasy to yield the invertibility \cite{vasyessential}. In their approach, the wave operator acts on \emph{anisotropic Sobolev spaces} (with order varying in phase space), and the Fredholm property is obtained by combining propagation of singularities estimates with \emph{radial estimates} at infinity (building on earlier works by Vasy \cite{kerrds} and Baskin--Vasy--Wunsch \cite{BVW}). The Feynman problem is then defined by imposing decrease of the regularity along the bicharacteristic flow, and it was shown by Vasy--Wrochna that this can in fact be  interpreted as a condition on data at past and future infinity in the sense of a \emph{geometric scattering theory} \cite{VW}. In this picture,  the splitting of the radial sets at infinity into \emph{sinks} and \emph{sources} plays the r\^{o}le of the asymptotic decomposition into particles and anti-particles.  Similar arguments have also been applied to asymptotically de Sitter spaces (the Lorentzian analogue of asymptotically hyperbolic spaces) \cite{GHV,positive,VW,bunchdavies}.

The Gell-Redman--Haber--Vasy method also allows for solving a \emph{non-linear Feynman problem} \cite{GHV}. Closely related techniques were recently used by Hassell--Gell-Redman--Schapiro--Zhang to prove the existence of standing waves for the non-linear Helmholtz equation \cite{hassell}, and it is expected that the Feynman problem for the Klein--Gordon equation could be solved by a mixture of the techniques in the two works \cite{GHV,hassell}. 
 
In the present notes we report on a different method which uses time-dependent pseudo-differential operators and techniques from scattering theory, applied to the case of the Klein--Gordon equation on asymptotically Minkowski space \cite{massivefeynman1,massivefeynman2}.  

\subsection{Setting and main results} On $M=\rr^{1+d}$, we consider an \emph{asymptotically Minkowski metric}  $g$ in the sense that:

\ben
\item[(1)] $g-g_{0}\in S^{-\delta}(\rr^{1+d})$ for some  $\delta>1$,
\item[(2)] $(M,g)$ is  {non-trapping},
\item[(3)] $(M,g)$ admits a time function  (i.e.~a smooth function with time-like gradient) that differs from the Minkowski $t$ coordinate by a term in $S^{1-\epsilon}(\rr^{1+d})$, $\epsilon>0$.
\een
Above, $g_{0}=-dt^2+dx^2$ is the Minkowski metric, and $ S^{-\delta}(\rr^{1+d})$ consists of metrics with coefficients behaving as a symbol of order $-\delta$ (thus, they decay in \emph{all} space-time directions at a $-\delta$ rate, and taking derivatives yields stronger decay).  The \emph{non-trapping} is the property that all null geodesics escape to infinity as the affine parameter tends to $+\infty$ or $-\infty$.

We consider the linear Klein--Gordon operator $-\Box_g + m^2$ (where $m>0$). In the first step, by constructing new coordinates and by composing $-\Box_g + m^2$ with multiplication operators, we show that the situation can be reduced to an operator $P$ on $\rr^{1+d}$ of the form:
\beq\label{main1}
P=\p_t^2 + A(t),
\eeq
 where $\rr\ni t \mapsto A(t)$ is a family of  differential operators on $\rr^d$ such that
\beq\label{main2}
A(t) -  \left(-\Delta_x + m^2\right) \in\Psi_{\rm sc}^{2,-\delta}(\rr; \rr^d).
\eeq
Here, $\Psi_{\rm sc}^{k,\ell}(\rr;\rr^d)$ is the space of $t$-dependent pseudo-differential operators obtained by quantizing $t$-dependent symbols $a(t,x,\xi)$ satisfying:
\[
\p_t^\gamma \p_x^\alpha \p_\xi^\beta a(t,x,\xi) \in O\big( (\bra t\ket+\bra x\ket)^{\ell-\gamma-|\alpha|} \bra \xi \ket^{k-|\beta|}\big).
\] 
for all $\gamma\in \nn_0$, $\alpha,\beta\in \nn_0^d$, where $\bra x \ket =(1+|x|^2)^\12$. Thus, \eqref{main2} means that $P$ differs from its Minkowski space analogue $P_{0}=\partial_t^2-\Delta_x + m^2$ by a term that decays simultaneously in $t$ and $x$ (the assumption $\delta>1$ means that this perturbation is \emph{short-range}). 

As already remarked, the Cauchy evolution $U_0(t,s)$ of $P_0$ is generated by the following Hamiltonian, denoted from now on by $H_0$:
\[
H_{0}= \begin{pmatrix} 0 & 1 \\ -\Delta_x + m^2 & 0 \end{pmatrix},
\]
which is self-adjoint in the energy space.  For some fixed $\gamma\in\open{\12,\12+ \delta}$ and  $m\in\rr$ we define:
\[
\cY =  \langle t\rangle^{-\gamma}L^{2}(\rr; H^{m}(\rr^{d})),
\]
and
\[
\cX = \big\{  u\in C^{0}(\rr; H^{m+1}(\rr^{d}))\cap C^{1}(\rr; H^{m}(\rr^{d})) \, : \, Pu\in \cY\big\}.
\]
Note that since $\langle t\rangle^{-\gamma}L^{2}(\rr)\subset L^{1}(\rr)$ for $\gamma>\12$, $\cY$ is a natural choice of Hilbert space in the context of solving an inhomogeneous Cauchy problem with $Pu\in \cY$, and then the solution $u$ lies in $\cX$.  We define the subspace of functions that satisfy \emph{Feynman conditions} at infinity:
\[
\cX_{\rm F} = \left\{  u\in \cX : \ \lim_{t\to\pm\infty}\one_{\rr^\mp}(H_0)\begin{pmatrix} u(t) \\ D_t u(t)\end{pmatrix} =0    \right\}.
\]
Our main result is the following theorem.
\begin{theorem}[\cite{massivefeynman1,massivefeynman2}] Let $P$ be as in \eqref{main1}--\eqref{main2} with $\delta>1$. Then $P: \cX_{\rm F}\to \cY$ is {invertible}. Its inverse $P_{\rm F}^{-1}$  coincides with the Duistermaat--H\"ormander  parametrix modulo smoothing terms. 
\end{theorem}

The theorem also applies to the original Klein--Gordon operator $-\Box_g+m^2$ on asymptotically Minkowski spacetimes, with minor  modifications to account for the change of coordinates.

\section{Sketch of the proof}

\subsection{Proof of Fredholm property} We start by showing the Fredholm property of  $P: \cX_{\rm F}\to \cY$. The main ingredient of the proof is the \emph{approximate diagonalisation} of the Cauchy evolution $U(t,s)$ of $P$.

In the first step, we use a new variant of the method used before in \cite{junker,GW,GOW,inout} to construct a family of elliptic pseudo-differential operators $\rr\ni t \mapsto B(t)$ such that:
\beq\label{eq:fac0}
B(t)=\sqrt{-\Delta_x+m^2} \mbox{ mod } \Psi_{\rm sc}^{1,-\delta}(\rr;\rr^d)
\eeq
and
\beq\label{eq:fac}
\bea  
P &= (D_t-B(t)) (D_t + B(t))  \mbox{ mod } \Psi_{\rm sc}^{-\infty,-1-\delta}(\rr;\rr^d)\\
 &= (D_t+B^*(t)) (D_t - B^*(t))  \mbox{ mod } \Psi_{\rm sc}^{-\infty,-1-\delta}(\rr;\rr^d).
\eea
\eeq
The terms in $\Psi_{\rm sc}^{-\infty,-1-\delta}(\rr;\rr^d)$ are smoothing and decaying, so they are merely compact errors that will not affect the Fredholm property (even though this is quite delicate considering that we mean compactness in the sense of the spaces $\cX_{\rm F}$ and $\cY$). The proof of \eqref{eq:fac} uses the $\Psi_{\rm sc}^{k,\ell}(\rr;\rr^d)$  symbolic calculus  and proceeds recursively, starting from the highest order terms in the poly-homogeneous expansion of $B(t)$.

Next, we show that it is possible to find $B(t)$ as above with the extra property that $(B+B^*)^{-1}(t)$ exists. This implies that if we set
\[
\bea
\widetilde{u}_1 & = (D_t+B(t))u,   \\
\widetilde{u}_2 & = (D_t - B^*(t)) u,
\eea
\]
then the transformation $u\mapsto \widetilde u=(\widetilde{u}_1 ,\widetilde{u}_2 )$ is invertible. By \eqref{eq:fac}, the equation $Pu=0$ is equivalent  \emph{modulo smoothing, decaying errors} to the diagonal system $\widetilde{P}\widetilde{u}=0$, where
\beq\label{diago}
\widetilde{P} = D_t - \widetilde{H}(t), \ \ \widetilde{H}(t)=\begin{pmatrix} B(t) & 0 \\ 0 & -B^*(t) \end{pmatrix}.
\eeq
Let us briefly comment on the invertibility of $(B+B^*)(t)$: it is arranged by modifying the original definition of $B(t)$ at low frequencies using the functional calculus. An important technical part of the argument consists then in proving that this spectral-theoretic operation is compatible with the time-dependent pseudo-differential calculus $\Psi_{\rm sc}^{k,\ell}(\rr;\rr^d)$ and does not affect  the two properties \eqref{eq:fac0}--\eqref{eq:fac}.

From now on we focus on the diagonal system \eqref{diago} (for the sake of brevity we disregard here a further transformation made in \cite{inout,massivefeynman1} which yields simpler preserved quantities). The analogue of the two function spaces $\cX$ and $\cY$ are the spaces
\[
\bea
\widetilde\cY &=  \langle t\rangle^{-\gamma}L^{2}(\rr; H^{m}(\rr^{d})\oplus H^{m}(\rr^{d})),\\
\widetilde\cX &= \big\{  \widetilde{u}\in C^{0}(\rr; H^{m}(\rr^{d})\oplus H^{m}(\rr^{d})): \ \widetilde{P}\widetilde{u}\in \widetilde\cY\big\}.
\eea
\]
The r\^{o}le of the two spectral projections $\one_{\rr^\pm}(H_0)$ is now simply played by the  projections $\pi^\pm$ to the respective two components of $\widetilde u=(\widetilde{u}_1 ,\widetilde{u}_2 )$. Therefore, the analogue of the space $\cX_{\rm F}$ is
\beq\label{eqfe}
\widetilde{\cX}_{\rm F} = \big\{  \widetilde{u}\in \widetilde\cX : \ \lim_{t\to\pm\infty} \pi^\mp \widetilde{u}(t)   =0    \big\}.
\eeq
This way, the Fredholm property of $P:  \cX_{\rm F}\to \cY$ is reduced  to the Fredholm property of the equivalent (modulo compact errors) problem $\widetilde{P}: \widetilde\cX_{\rm F}\to \widetilde\cY$. Note that dealing with the compact errors is not completely straightforward because $P$, its approximately diagonalized version and the diagonal operator $\widetilde P$ act on different spaces;  for the sake of brevity we omit the details here.

\begin{proposition}\label{prop:fredholm} The operator $\widetilde{P}: \widetilde\cX_{\rm F}\to \widetilde\cY$ is Fredholm.
\end{proposition}

In fact, we can even construct an explicit \emph{inverse} of $\widetilde{P}$. We set 
\beq\label{eq:tildeinverse}
\big(\widetilde{P}^{-1}_{\rm F} \widetilde{u}\big)(t)  =\int_{-\infty}^t \pi^+ \widetilde{U}(t,s)\widetilde{u}(s) ds- \int_{t}^{\infty} \pi^- \widetilde{U}(t,s)\widetilde{u}(s) ds,
\eeq
where $\widetilde{U}(t,s)$ is the Cauchy evolution of $\widetilde{P}$ or equivalently, the evolution generated by the time-dependent Hamiltonian $\widetilde{H}(t)$.

 To show that $\widetilde{P}^{-1}_{\rm F}$ is a bounded operator from $\widetilde{\cY}$ to $\widetilde{\cX}_{\rm F}$, we have to check that the RHS of \eqref{eq:tildeinverse} satisfies the  Feynman conditions at infinity, and in particular we need a good control of  ${t\to\pm \infty}$ limits. To that end, observe that the time-dependent Hamiltonian $\widetilde{H}(t)$ is asymptotic to the operator
 $$
\widetilde{H}_0=\begin{pmatrix} \sqrt{-\Delta_x + m^2} & 0 \\ 0 & -\sqrt{-\Delta_x + m^2} \end{pmatrix}, 
 $$
 and so the evolution $\widetilde{U}_0(t,s)=e^{i(t-s)\widetilde{H}_0}$ generated by $\widetilde H_0$ provides a natural comparison dynamics. Now, the idea is that because $\pi^\pm$ commutes with all our diagonal operators, the Feynman conditions $\lim_{t\to\pm\infty} \pi^\mp v(t)   =0$ are equivalent to 
 $$
 \lim_{t\to\pm\infty} \pi^\mp \widetilde{U}_0(s,t) v(t)   =0 .
 $$ 
Thus, checking that \eqref{eq:tildeinverse} satisfies these conditions amounts to a good control of $\widetilde{U}_0(s,t) \widetilde{U}(t,s)$ at large $|t|$. But instead of doing that directly, we can deduce it from the behaviour of the $t$ derivative: this produces an $\widetilde{H}(t)-\widetilde{H}_0$ factor, and the good news is that the time-decay of this difference is  well under control thanks to \eqref{eq:fac0}. 

Finally, once we know that  $\widetilde{P}^{-1}_{\rm F}:\widetilde{\cY}\to\widetilde{\cX}_{\rm F}$ it is straightforward to check that it is an inverse

We refer to  \cite{massivefeynman1} for the detailed proof.
 

\subsection{Proof of invertibility} We have proved the Fredholm property of  $P:\cX_{\rm F}\to \cY$ by showing its invertibility modulo compact terms (which follows from  invertibility of $\widetilde{P}:\widetilde{\cX}_{\rm F}\to \widetilde{\cY}$), and so in particular its Fredholm index is $0$. Therefore, to conclude that it is invertible we only need to prove that $\Ker P=\{0\} $. 

Let us first illustrate the argument in the diagonalized setting (even though in that case the invertibility can be seen more directly).  Let $\widetilde u \in \Ker \widetilde P$, understood as a subspace of $\widetilde{\cX}_{\rm F}$. On the one hand side, we show that this implies
\beq\label{e1.3}
\begin{array}{rl}
	&\| \pi^{\pm}\widetilde u(t)\|^{2}_{L^2\oplus L^2} = O(t^{1- \delta}) \hbox{ as }t\to \mp \infty,\\[2mm]
	&\| \pi^{\pm} \widetilde u(t)\|^{2}_{L^2\oplus L^2}= c^{\pm}+ O(t^{1-\delta})\hbox{ as  }t\to \pm \infty
\end{array}
\eeq
for some constants $c^\pm$. This uses arguments from scattering theory similar to the proof of Proposition \ref{prop:fredholm}. On the other hand, following an argument due to Vasy \cite{vasyessential}, for a well-chosen  family $\{\chi_\epsilon\}_{\epsilon>0}$ of smooth functions in $t$ satisfying $\supp\chi_\epsilon\subset \{|t|\leq 2\epsilon^{-1}\}$, by integration by parts we find:
\beq\label{eq:brackets}
\bea
0&= \bra \widetilde P \widetilde u, \chi_\epsilon \widetilde u \ket  - \bra \chi_\epsilon \widetilde u , \widetilde P \widetilde u\ket \\
&= - i \bra \widetilde u,\p_t \chi_\epsilon \widetilde u \ket
\eea
\eeq
for a suitable inner product for which $\widetilde P$ is formally self-adjoint.
We then split the RHS into a term with $\pi^+$ and a term with $\pi^-$. By positivity properties of $\pi^\pm$ and taking into account the $\epsilon$ behaviour of $\p_t \chi_\epsilon$, we deduce from \eqref{eq:brackets}  that the constants $c^\pm$ in \eqref{e1.3} are of order $O(\epsilon^{\delta-1})$. Since $\epsilon$ can be taken arbitrarily close to $0$, it follows that $\lim_{t\to\infty}\widetilde u(t)=0$ in $L^2(\rr^d)\oplus L^2(\rr^d )$. Therefore, $\widetilde u$ is a solution with vanishing scattering data, and consequently $\widetilde u =0$. 

Crucially, this conclusion is also  valid for the \emph{approximately} diagonalized operator to which $P$ is equivalent. Indeed, \emph{modulo smoothing, decaying terms},  the approximately diagonal operator equals $\widetilde{P}$ (acting on slightly different spaces, to which the argument above still applies). Therefore, the asymptotic properties of  solutions are the same as for $\widetilde{P}$, and a perturbative argument can be applied.  This implies  that $\Ker P=\{0\} $ as asserted. 

The  detailed proof can be found in  \cite{massivefeynman2}.

\subsection{Proof of microlocal properties} Using Fourier Integral Operators techniques, Duistermaat and H\"{o}rmander have constructed a Feynman parametrix $H_{\rm F}$ in the sense of being an inverse of $P$ (on say, test functions) modulo smoothing terms (but not necessarily compact) \cite{DH}. They have shown that $H_{\rm F}$ is characterized  \emph{uniquely modulo smoothing term} by the \emph{wavefront set} of its Schwartz kernel. Therefore, one way of proving that our inverse $P_{\rm F}^{-1}$ coincides with $H_{\rm F}$ modulo smoothing terms is to show that the wavefront sets are equal. 

Let us recall  that if $u$ is a distribution, its wavefront set $\wf(u)\subset T^*M$ is defined as follows:  $(t,x,\tau,\xi)\in T^*M\setminus \zero$ is \emph{not} in  $\wf(u)$ if there exists a pseudo-differential operator $A$ (say, properly supported) such that its principal symbol $\sigma_{\rm pr}(A)$ is non-zero at $(t,x,\tau,\xi)$ and $Au\in \cf(M)$. In particular, if $Pu=0$ then $\wf(u)\subset \sigma_{\rm pr }(P)^{-1}(\{0\})$ (the \emph{characteristic set} of $P$). In our case of interest the characteristic set has two connected components, denoted in what follows by $\cN^+$ and $\cN^-$. For instance, on Minkowski space, $\cN^\pm=\{ \tau = \pm |x|\}$ and it is not difficult to see that solutions of the half-Klein--Gordon equation $(D_t \mp \sqrt{-\Delta_x + m^2} )u=0$ propagate with wavefront set in $\cN^\pm$.

We only very briefly explain the idea  of the proof. Without writing the explicit formula, we stress that what is specific about the wavefront set of $H_{\rm F}$ is that it \emph{distinguishes between the direction of propagation in the two components} $\cN^\pm$ of the characteristic set of $P$.

First of all, the transformation that, modulo smoothing terms,  replaces $P$ by  $\widetilde{P}$,  does not affect wave front sets. In consequence we are reduced to estimating the wavefront set of  $\widetilde{P}_{\rm F}^{-1}$, defined in formula \eqref{eq:tildeinverse}. The crucial observation is that in that formula, disregarding components that are identically zero, $\pi^+ \widetilde U(t,s)$ is  the  evolution generated by the elliptic family $B(t)$ with \emph{positive} principal symbol, whereas $\pi^- \widetilde U(t,s)$ is generated by the  family $-B^*(t)$ with \emph{negative} principal symbol. This corresponds precisely to propagation within $\cN^\pm$, and the time integrals in \eqref{eq:tildeinverse} account for the correct direction of propagation. This information can be then used to deduce the full wavefront set of $\widetilde{P}_{\rm F}^{-1}$ and thus of $P_{\rm F}^{-1}$.

The detailed proof can be found in \cite{massivefeynman1}.

\section{Further related results}

\subsection{Related recent works and outlook}  We conclude the discussion by mentioning a couple of further related results.

We remark that as already indicated in the introduction,  the Feynman problem  is closely related to the limiting absorption principle for $P=-\Box_g +m^2$. The essential self-adjointess of $P$ and the fact that $P_{\rm F}^{-1}=\lim_{\varepsilon\to 0^+}(P-i \varepsilon)^{-1}$ in the sense of operators acting between suitable Hilbert spaces (different from $L^2(M,g)$) was shown by Derezi\'nski--Siemssen in the  case of \emph{static} spacetimes \cite{derezinskisiemssen0,derezinskisiemssen1}. The more difficult case of \emph{non-trapping Lorentzian scattering spaces} (closely related to the asymptotically Minkowski spacetimes considered here)  was solved by Vasy \cite{vasyessential} (improving this way an earlier  result on positivity of propagator differences \cite{positive}). Very recently, Nakamura--Taira gave a proof of essential-self-adjointness for a broad class of metrics with asymptotically constant coefficients \cite{nakamura}. We also mention earlier works on Strichartz estimates  \cite{tzvetkov,taira}, which  include essential self-adjointness results under more restrictive assumptions. 

One of the motivations for the Feynman problem comes from its similarity to a result of B\"{a}r--Strohmaier, who proved the Fredholm property of the massless Dirac operator with  conditions analogous to the ones considered here, but at \emph{finite} times \cite{BS}. These can be then interpreted as a Lorentzian version of the Atiyah-Patodi-Singer boundary conditions, and the work \cite{BS} provides indeed a geometric formula for the index, which can also be interpreted in Quantum Field Theory as the created charge \cite{BS2}. 

A more abstract version of our setting was recently considered by Derezi\'nski--Siemssen \cite{derezinskisiemssen1}, who gave sufficient conditions for Feynman inverses, formulated in terms of the asymptotic spectral projections and their image by scattering operators. The work  \cite{derezinskisiemssen1}  also discusses Fock space implentability  and explains  how the Feynman propagator provides a rigorous interpretation of the \emph{in-out expectation value of time-ordered products of fields} often considered by physicists.

In summary, this collection of results demonstrates that the Feynman propagator or inverse is the closest Lorentzian analogue of the boundary value of the Laplacian resolvent. Its important r\^{o}le in quantum physics, its ties to the global geometry and its possible relationships with local geometric quantities makes it an extremely interesting object for further investigations.

 \end{document}